\documentclass{amsproc}

\newtheorem{theorem}{Theorem}[section]
\newtheorem{corollary}[theorem]{Corollary}
\newtheorem{proposition}[theorem]{Proposition}

\newcommand{\I}{\mathbb{I}}
\newcommand{\N}{\mathbb{N}}
\newcommand{\T}{\mathbb{T}}
\newcommand{\Z}{\mathbb{Z}}
\newcommand{\cW}{\mathcal{W}}
\begin{document}

\title[Algebras of Functions with Fourier Coefficients in Orlicz Spaces]
{Algebras of Functions with {F}ourier Coefficients \\
in Weighted {O}rlicz Sequence Spaces}

\author[Alexei~Yu.~Karlovich]{{\bf Alexei Yu. Karlovich}}

\address{
Alexei~Yu.~Karlovich\\
Departamento de Mathem\'atica\\
Instituto Superior T\'ecnico\\
Av. Rovisco Pais\\
1049-001 Lisboa \\
Portugal}
\email{akarlov@math.ist.utl.pt}

\thanks{The author is supported by F.C.T. (Portugal)
grant SFRH/BPD/11619/2002.}

\dedicatory{To the Memory of Erhard Meister}

\begin{abstract}
We prove that the set of all integrable functions whose
sequences of negative (resp. nonnegative) Fourier coefficients belong to 
$\ell^1\cap\ell^\Phi_{\varphi,w}$ (resp. to $\ell^1\cap\ell^\Psi_{\psi,\varrho}$), 
where $\ell^\Phi_{\varphi,w}$ and $\ell^\Psi_{\psi,\varrho}$ are
two-weighted Orlicz sequence spaces, forms an algebra under pointwise 
multiplication whenever the weight sequences 
\[
\varphi=\{\varphi_n\},\quad
\psi=\{\psi_n\},\quad
w=\{w_n\},\quad
\varrho=\{\varrho_n\} 
\]
increase and satisfy the $\Delta_2$-condition. 
\end{abstract}

\maketitle
\section{Introduction}
Let $\T$ be the unit circle. For a complex-valued function $f\in L^1(\T)$, let
$\{f_n\}_{n\in\Z}$ be the sequence of the Fourier coefficients of $f$,
\[
f_n:=\frac{1}{2\pi}\int_0^{2\pi}f(e^{i\theta})e^{-in\theta}d\theta.
\]
Let $W$ be the Wiener algebra of all functions $f$ on $\T$ for which
\[
\|f\|_W:=\sum_{k=-\infty}^\infty |f_k|<\infty.
\]
Let $F\ell^{p,r}_{\alpha,\beta}$, where
$1\le p,r<\infty$ and $0\le\alpha,\beta<\infty$,
denote the set of all functions $f$ on $\T$ for which
\[
\|f\|_{F\ell^{p,r}_{\alpha,\beta}}:=
\left(\sum_{k=1}^\infty |f_{-k}|^p(k+1)^{\alpha p}\right)^{1/p} 
+
\left(\sum_{k=0}^\infty |f_k|^r(k+1)^{\beta r}\right)^{1/r}<\infty.
\]
The following result was conjectured by A.~B\"ottcher and 
B.~Silbermann and proved by Detlef Horbach.
The proof is in \cite[Section~6.54]{BS90}.
\begin{theorem}\label{th:Horbach}
If $1\le p,r<\infty$ and $0\le\alpha,\beta<\infty$, then 
$W\cap F\ell^{p,r}_{\alpha,\beta}$ is an algebra under pointwise 
multiplication.
\end{theorem}
Note that a stronger fact is actually proved in \cite[Section~6.54]{BS90}:
the inequality 
$\|f^2\|_{F\ell^{p,r}_{\alpha,\beta}}\le C_1\|f\|_W\|f\|_{F\ell^{p,r}_{\alpha,\beta}}$
is fulfilled for some $C_1>0$ and any $f\in W\cap F\ell^{p,r}_{\alpha,\beta}$. From the latter
inequality and $\|fg\|_W\le \|f\|_W\|g\|_W$ one can easily get
\[
\|fg\|_W+\|fg\|_{F\ell^{p,r}_{\alpha,\beta}}
\le C
\Big(\|f\|_W+\|f\|_{F\ell^{p,r}_{\alpha,\beta}}\Big)
\Big(\|g\|_W+\|g\|_{F\ell^{p,r}_{\alpha,\beta}}\Big)
\]
for some $C>0$ and all $f,g\in W\cap F\ell^{p,r}_{\alpha,\beta}$,
that is, $W\cap F\ell^{p,r}_{\alpha,\beta}$ is a Banach algebra.
This result has important applications in 
the theory of Toeplitz determinants and
the theory of Toeplitz operators on Lebesgue spaces 
$l^p_\varrho (1<p<\infty)$ with Khvedelidze weights 
$\varrho$ (see \cite[Ch.~7 and Ch.~10]{BS90}).
The aim of the present note is to generalize this result  
by replacing in the definition of  $F\ell^{p,r}_{\alpha,\beta}$
\begin{enumerate}
\item
norms in weighted Lebesgue sequence
spaces by norms in two-weighted Orlicz sequence spaces;
\item
canonical weight sequences $\{(n+1)^\lambda\}_{n=0}^\infty$ 
with $\lambda\in\{\alpha,\beta\}$ by general increasing weight 
sequences satisfying the $\Delta_2$-condition.
\end{enumerate}
The author hopes that such an extension will be useful because the
scale of Orlicz spaces is much wider and sensitive than the scale
of Lebesgue spaces.

Note also that similar questions were considered by P.~L.~Ul'yanov 
\cite{Ulyanov92} for the set $E$ of functions $f\in C(\T)$ with
\[
\|f\|_E:=\sum_{k\in\Z}\omega(|f_k|)\tau_n<\infty,
\]
where $\omega$ is a function with certain properties (a so-called modulus 
of continuity) and $\{\tau_n\}_{n\in\Z}$ is a weight sequence. Under some 
natural assumptions he proved that $E$ is an algebra under pointwise 
multiplication.

This note is organized as follows. In Section~\ref{sec:2}
we remind the notion of two-weighted Orlicz sequence spaces and define
a related class $F\ell^{\Phi,\Psi}_{\varphi,w;\psi,\varrho}$ that generalizes
the class $F\ell^{p,r}_{\alpha,\beta}$. Further we formulate the main result:
the set $W\cap F\ell^{\Phi,\Psi}_{\varphi,w;\psi,\varrho}$ is a Banach algebra 
under pointwise multiplication whenever the weight sequences 
$\varphi,w,\psi,\varrho$ increase and satisfy the $\Delta_2$-condition. 
Finally we state a corollary about factorization 
of nondegenerate functions in this algebra. Section~\ref{sec:3} 
contains the proofs.
\section{Preliminaries and the main result}\label{sec:2}
\subsection{Weighted Orlicz sequence spaces}
An \textit{Orlicz function} $\Phi$ is a continuous non-decreasing
and convex function defined for $t\ge 0$ such that $\Phi(0)=0$
and $\lim\limits_{t\to\infty}\Phi(t)=\infty$ (see \cite[Definition~4.a.1]{LT77}).
Let $\I$ be either $\N:=\{1,2,\dots\}$ or $\Z_+:=\N\cup\{0\}$ and let $\{\Phi_n\}_{n\in\I}$
be a sequence of Orlicz functions. The set $\ell^{\{\Phi_n\}}(\I)$ of all
sequences of complex numbers $c=\{c_n\}_{n\in\I}$ with
\[
\sum_{n\in\I}\Phi_n\left(\frac{|c_n|}{\lambda}\right)<\infty
\]
for some $\lambda=\lambda(c)>0$ is a Banach space when equipped with the norm
\[
\|c\|_{\ell^{\{\Phi_n\}}(\I)}
=
\inf\left\{\lambda>0:\quad 
\sum_{n\in\I}\Phi_n\left(\frac{|c_n|}{\lambda}\right)\le 1
\right\}.
\]
The space $\ell^{\{\Phi_n\}}(\I)$ is called a \textit{modular sequence space}
(see \cite[Definition~4.d.1]{LT77}) or a \textit{Musielak-Orlicz sequence 
space} (see, e.g., \cite{Musielak83}).

Any sequence $\{\nu_n\}_{n\in\I}$ of positive numbers is called a 
\textit{weight sequence}. Let $\Phi$ be an Orlicz function and 
$\varphi=\{\varphi_n\}_{n\in\I}$ and $w=\{w_n\}_{n\in\I}$ be weight 
sequences. Clearly 
\begin{equation}\label{eq:Orlicz-sequence}
\Phi_n(x):=\Phi(x\varphi_n)w_n \quad (x\ge 0)
\end{equation}
is an Orlicz function for every $n\in\I$. 
The special case of a Musielak-Orlicz sequence space generated by the 
sequence of Orlicz functions (\ref{eq:Orlicz-sequence}) is denoted by
$\ell^\Phi_{\varphi,w}(\I)$ and is called the \textit{two-weighted Orlicz space}
generated by the Orlicz function $\Phi$ and the weight sequences 
$\varphi=\{\varphi_n\}_{n\in\I}$ and $w=\{w_n\}_{n\in\I}$.

In particular, if $\Phi(x)=x^p, 1\le p<\infty$, then $\ell^\Phi_{\varphi,w}(\I)$
is the weighted Lebesgue sequence space $\ell^p_\rho(\I)$ and
\[
\|c\|_{\ell^\Phi_{\varphi,w}(\I)}=\left(\sum_{n\in\I} |c_n|^p \rho_n\right)^{1/p}
\quad\quad
\Big(\rho_n=\varphi_n^p w_n,\quad n\in\I\Big).
\]
So, in this case, $\ell^\Phi_{\varphi,1}(\I)=\ell^\Phi_{1,\varphi^p}(\I)$, where
$\varphi^p:=\{\varphi_n^p\}_{n\in\I}$.
But in the general case an Orlicz function $\Phi$ is not homogeneous. Hence
the weighted Orlicz sequence spaces $\ell^\Phi_{\varphi,1}(\I)$
and $\ell^\Phi_{1,w}(\I)$ are essentially different.

Applying criteria of coincidence of Musielak-Orlicz sequence spaces
and Musi\-elak-Orlicz sequence classes (see \cite[Theorem~8.13(b)]{Musielak83})
to the case of the sequence of Orlicz functions (\ref{eq:Orlicz-sequence}),
one can get the following.
\begin{proposition}
Suppose $\Phi$ is an Orlicz function, $\Phi(x)>0$ whenever $x>0$, and
$\varphi=\{\varphi_n\}_{n\in\I},w=\{w_n\}_{n\in\I}$ are weight
sequences. Then
\[
\ell^\Phi_{\varphi,w}(\I)=\left\{
c=\{c_n\}_{n\in\I}:\quad
\sum_{n\in\I}\Phi(|c_n|\varphi_n)w_n<\infty
\right\}
\]
if and only if there exist positive numbers $\delta,K$ and a sequence of 
nonnegative numbers $\{d_n\}_{n\in\I}$ such that for all $x\ge 0$ and
$n\in\I$,
\[
\Phi(x\varphi_n)w_n<\delta
\quad\Longrightarrow\quad
\Phi(2x\varphi_n)w_n\le K\Phi(x\varphi_n)w_n+d_n
\]
and $\sum_{n\in\I} d_n<\infty$.
\end{proposition}
\subsection{The main result}
Let $\Phi$ and $\Psi$ be Orlicz functions and let 
\[
\varphi=\{\varphi_n\}_{n=1}^\infty,\quad
w=\{w_n\}_{n=1}^\infty,\quad
\psi=\{\psi_n\}_{n=0}^\infty,\quad
\varrho=\{\varrho_n\}_{n=0}^\infty
\]
be weight sequences. We denote by $F\ell^{\Phi,\Psi}_{\varphi,w;\psi,\varrho}$
the set of all functions $f\in L^1(\T)$ such that the sequence 
$\{f_{-n}\}_{n\in\N}$ of all negative Fourier coefficients of $f$ belongs
to the two-weighted Orlicz space $\ell^\Phi_{\varphi,w}(\N)$
and the sequence $\{f_n\}_{n\in\Z_+}$ of all nonnegative Fourier coefficients
of $f$ belongs to the two-weighted Orlicz space $\ell^\Psi_{\psi,\varrho}(\Z_+)$.
The set $W\cap F\ell^{\Phi,\Psi}_{\varphi,w;\psi,\varrho}$ is a Banach space
with respect to the norm
\[
\|f\|_{W\cap F}
:=
\|f\|_W+\|f\|_-+\|f\|_+,
\]
where
\begin{eqnarray*}
\|f\|_-
&:=&
\Big\|\{f_{-k}\}_{k\in\N}\Big\|_{\ell^\Phi_{\varphi,w}(\N)}
=
\inf\left\{\lambda>0:\quad
\sum_{k=1}^\infty \Phi\left(\frac{|f_{-k}|\varphi_k}{\lambda}\right)w_k
\le 1\right\},
\\
\|f\|_+
&:=&
\Big\|\{f_k\}_{k\in\Z_+}\Big\|_{\ell^\Psi_{\psi,\varrho}(\Z_+)}
=
\inf\left\{\mu>0:\quad
\sum_{k=0}^\infty \Psi\left(\frac{|f_k|\psi_k}{\mu}\right)\varrho_k\le 1\right\}.
\end{eqnarray*}

We denote by $\cW^+$ (resp. by $\cW^-$) the collection of all weight sequences
$\{\nu_n\}_{n=0}^\infty$ (resp. $\{\nu_n\}_{n=1}^\infty$) such that
\begin{itemize}
\item[(i)]
$\nu_0>0$ (resp. $\nu_1>0$);
\item[(ii)]
$\nu_n\le\nu_{n+1}$ for $n\in\Z_+$ (resp. for $n\in\N$);
\item[(iii)]
$\{\nu_n\}_{n=0}^\infty$ (resp. $\{\nu_n\}_{n=1}^\infty$) satisfies
the $\Delta_2$-condition, that is, there exists a constant $C_\nu\in(0,\infty)$
such that $\nu_{2n}\le C_\nu\nu_n$ for $n\in\N$.
\end{itemize}
From (ii) and (iii) it follows that $C_\nu\ge 1$.
\begin{theorem}\label{th:main}
If $\Phi,\Psi$ are arbitrary Orlicz functions, 
$\varphi=\{\varphi_n\}_{n=1}^\infty, w=\{w_n\}_{n=1}^\infty$
are weight sequences in  $\cW^-$, 
and $\psi=\{\psi_n\}_{n=0}^\infty,\varrho=\{\varrho_n\}_{n=0}^\infty$
are weight sequences in $\cW^+$, then for every 
$f,g\in W\cap F\ell^{\Phi,\Psi}_{\varphi,w;\psi,\varrho}$,
\begin{equation}\label{eq:main*}
\|fg\|_{W\cap F}\le C\|f\|_{W\cap F}\|g\|_{W\cap F}
\end{equation}
where $C:=1+2(1+C_w)C_\varphi+2(1+C_\varrho)C_\psi$.
\end{theorem}

This theorem will be proved in the next section.

Clearly, the weight sequences $\widetilde{\varphi}=\{(n+1)^\alpha\}_{n=1}^\infty$ and 
$\widetilde{\psi}=\{(n+1)^\beta\}_{n=0}^\infty$ belong to $\cW^-$ and $\cW^+$, 
respectively, whenever $0\le\alpha,\beta<\infty$. If 
\[
\widetilde{\Phi}(x):=x^p,
\quad
\widetilde{\Psi}(x):=x^r,
\quad
x\ge 0,
\quad
1\le p,r<\infty, 
\]
then $F\ell^{\widetilde{\Phi},\widetilde{\Psi}}_{\widetilde{\varphi},1;\widetilde{\psi},1}$
with the norm $\|f\|_F:=\|f\|_-+\|f\|_+$
is isometrically isomorphic to $F\ell^{p,r}_{\alpha,\beta}$. Hence 
Theorem~\ref{th:Horbach} follows from Theorem~\ref{th:main}.
\begin{corollary}
Under the assumptions of Theorem~{\rm \ref{th:main}} we have the following.

{\rm (a)}
$W\cap F\ell^{\Phi,\Psi}_{\varphi,w;\psi,\varrho}$ is a commutative
Banach algebra under pointwise multiplication.
The maximal ideal space of this algebra coincides with $\T$.

{\rm (b)}
If $b\in W\cap F\ell^{\Phi,\Psi}_{\varphi,w;\psi,\varrho}$ does not
vanish on $\T$ and the Cauchy index of $b$ vanishes, then $b$ has 
a logarithm in $W\cap F\ell^{\Phi,\Psi}_{\varphi,w;\psi,\varrho}$.

{\rm (c)}
If we let 
\[
G(b):=\exp\Big((\log b)_0\Big),
\quad\quad
b_\pm(t):=\exp\left(
\sum_{n=1}^\infty(\log b)_{\pm n}t^{\pm n}
\right)
\quad (t\in\T),
\]
then $b=G(b)b_-b_+$ and 
$b_\pm^{\pm 1}\in W\cap F\ell^{\Phi,\Psi}_{\varphi,w;\psi,\varrho}$.
\end{corollary}
\begin{proof}
By Theorem~\ref{th:main}, $W\cap F\ell^{\Phi,\Psi}_{\varphi,w;\psi,\varrho}$
is a commutative Banach algebra under pointwise multiplication.
The description of its maximal ideal space is standard.
Part (b) follows from part (a) and \cite[Section~2.41(e)]{BS90}.
Part (c) is an immediate consequence of part (b).
\end{proof}
\section{Proofs}\label{sec:3}
\subsection{Auxiliary results}
\begin{proposition}\label{pr:3}
If $\{\nu_n\}_{n=0}^\infty\in\cW^+$ (resp. $\{\nu_n\}_{n=1}^\infty\in\cW^-$)
and $k\in\Z_+$ (resp. $k\in\N$), then
\begin{equation}\label{eq:3.1}
\nu_k\le C_\nu\nu_j
\quad\mbox{for}\quad j\ge k-[k/2].
\end{equation}
\end{proposition}
\begin{proof}
If $k<2$, then $k-[k/2]=k$. Hence from (ii) and $C_\nu\ge 1$ we get 
(\ref{eq:3.1}).

If $k\ge 2$, then $2[k/2]\le k$. Therefore $k\le 2(k-[k/2])$. In that case
from (ii) and (iii) it follows that
\[
\nu_k\le\nu_{2(k-[k/2])}\le C_\nu\nu_{k-[k/2]}\le C_\nu\nu_j
\]
for $j\ge k-[k/2]$, i.e., we obtain (\ref{eq:3.1}).
\end{proof}

Let $\{f_n\}_{n=-\infty}^\infty$ and $\{g_n\}_{n=-\infty}^\infty$
be the Fourier coefficient sequences of functions $f\in L^1(\T)$  
and $g\in L^1(\T)$, respectively. Put $a_n:=|f_n|$ and $b_n:=|g_n|$.
\begin{proposition}\label{pr:1}
{\rm (a)} If $k\in\N$, then
\begin{equation}\label{eq:1.1}
|(fg)_{-k}|\le 
\sum_{j=0}^\infty a_jb_{-k-j}
+
\sum_{j=1}^{[k/2]}a_{-j}b_{-k+j}
+
\sum_{j=0}^\infty b_ja_{-k-j}
+
\sum_{j=1}^{[k/2]}b_{-j}a_{-k+j}.
\end{equation}
{\rm (b)} If $k\in\Z_+$, then
\[
|(fg)_k|\le 
\sum_{j=1}^\infty  a_{-j}b_{k+j}
+
\sum_{j=0}^{[k/2]} a_jb_{k-j}
+
\sum_{j=1}^\infty  b_{-j}a_{k+j}
+
\sum_{j=0}^{[k/2]} b_ja_{k-j}.
\]
\end{proposition}
\begin{proof}
(a) For $k\in\N$,
\begin{eqnarray}
\label{eq:1.2}
|(fg)_{-k}| 
&=& 
\left|\sum_{j=-\infty}^\infty f_jg_{-k-j}\right|
\le
\sum_{j=-\infty}^\infty a_jb_{-k-j}
\\
\nonumber
&=&
\sum_{j=-\infty}^{-k} a_jb_{-k-j}
+
\sum_{j=-k+1}^{-1} a_jb_{-k-j}
+
\sum_{j=0}^\infty a_jb_{-k-j}
\\
\nonumber
&=:&
\sigma_1+\sigma_2+\sigma_3.
\end{eqnarray}
Changing variables in $\sigma_1$ ($r=-k-j$), we obtain
\begin{equation}\label{eq:1.3}
\sigma_1:=\sum_{j=-\infty}^{-k}a_jb_{-k-j}
=
\sum_{r=0}^\infty a_{-k-r}b_r.
\end{equation}
Obviously, $\sigma_2=0$ if $k=1$. Hence (\ref{eq:1.1}) follows from
(\ref{eq:1.2}) and (\ref{eq:1.3}) for $k=1$.

If $k>1$ then changing variables in $\sigma_2$ ($r=-j$) we get
\[
\sigma_2:=\sum_{j=-k+1}^{-1} a_jb_{-k-j}=\sum_{r=1}^{k-1} a_{-r}b_{-k+r}.
\]

If $k=2m$ and $m\in\N$, then
\begin{equation}\label{eq:1.4}
\sigma_2
=
\sum_{j=1}^{2m-1} a_{-j}b_{-2m+j}
=
\sum_{j=1}^m a_{-j} b_{-2m+j}
+
\sum_{j=m+1}^{2m-1} a_{-j}b_{-2m+j}.
\end{equation}
Changing variables in the second sum ($r=2m-j$), we obtain
\begin{equation}\label{eq:1.5}
\sum_{j=m+1}^{2m-1} a_{-j}b_{-2m+j}
=
\sum_{r=1}^{m-1} a_{-2m+r}b_{-r}
\le
\sum_{j=1}^m b_{-j}a_{-2m+j}.
\end{equation}
Since $[k/2]=[2m/2]=m$, from (\ref{eq:1.2})--(\ref{eq:1.5}) we deduce 
that (\ref{eq:1.1}) holds for $k=2m$ and $m\in\N$.

If $k=2m+1$ and $m\in\N$, then
\begin{equation}\label{eq:1.6}
\sigma_2 = 
\sum_{j=1}^{2m} a_{-j}b_{-(2m+1)+j}
=
\sum_{j=1}^m a_{-j}b_{-(2m+1)+j}
+
\sum_{j=m+1}^{2m} a_{-j}b_{-(2m+1)+j}.
\end{equation}
Changing variables in the second sum ($r=2m+1-j$), we infer that
\begin{equation}\label{eq:1.7}
\sum_{j=m+1}^{2m} a_{-j}b_{-(2m+1)+j}
=
\sum_{r=1}^m a_{-(2m+1)+r}b_{-r}.
\end{equation}
Since $[k/2]=[(2m+1)/2]=m$, from (\ref{eq:1.2})--(\ref{eq:1.3}) and
(\ref{eq:1.6})--(\ref{eq:1.7}) we conclude that (\ref{eq:1.1}) is
satisfied for $k=2m+1$ and $m\in\N$. Part (a) is proved.
Part (b) is proved analogously to Part (a).
\end{proof}
\subsection{Proof of Theorem~\ref{th:main}}
\begin{proof}
The idea of this proof is borrowed from \cite[Theorem~6.54]{BS90}.
Let us show that
\begin{equation}\label{eq:main.0}
\|fg\|_-\le C_-(\|f\|_W\|g\|_-+\|g\|_W\|f\|_-),
\end{equation}
where $C_-:=(1+C_w)C_\varphi$, $C_\varphi$ and $C_w$ are the constants in 
the $\Delta_2$-condition for the sequences $\{\varphi_n\}_{n=1}^\infty$ 
and $\{w_n\}_{n=1}^\infty$, respectively. If
$\|f\|_W\|g\|_-+\|g\|_W\|f\|_-=0$, then (\ref{eq:main.0}) is obvious.

Assume that $\|f\|_W\|g\|_-+\|g\|_W\|f\|_->0$. Since $\Phi$ is 
increasing, from Proposition~\ref{pr:1}(a) it follows that for $k\in\N$, 
\begin{equation}\label{eq:main.1}
\Phi\left(
\frac{|(fg)_{-k}|\varphi_k}{C_-(\|f\|_W\|g\|_-+\|g\|_W\|f\|_-)}
\right)
\le
\Phi\left(
\sum_{i=1}^6\Sigma_i\Big/\sum_{i=1}^6\sigma_i
\right),
\end{equation}
where 
\[
\begin{array}{lcl}
\displaystyle
\Sigma_1 := 
\sum_{j=0}^\infty(a_j\|g\|_-)\frac{b_{-k-j}\varphi_k}{C_-\|g\|_-},
& &
\displaystyle
\sigma_1 :=
\sum_{j=0}^\infty a_j\|g\|_-, 
\\
\\
\displaystyle
\Sigma_2 :=
\sum_{j=1}^{[k/2]}(a_{-j}\|g\|_-)\frac{b_{-k+j}\varphi_k}{C_-\|g\|_-},
& &
\displaystyle
\sigma_2 :=
\sum_{j=1}^{[k/2]} a_{-j}\|g\|_-, 
\\
\\
\displaystyle
\Sigma_3 :=
\sum_{j<-[k/2]}(a_j\|g\|_-)\cdot 0=0,
& &
\displaystyle
\sigma_3 :=
\sum_{j<-[k/2]} a_j\|g\|_-,
\\
\\
\displaystyle
\Sigma_4 := 
\sum_{j=0}^\infty(b_j\|f\|_-)\frac{a_{-k-j}\varphi_k}{C_-\|f\|_-},
& &
\displaystyle
\sigma_4 :=
\sum_{j=0}^\infty b_j\|f\|_-, 
\\
\\
\displaystyle
\Sigma_5 :=
\sum_{j=1}^{[k/2]}(b_{-j}\|f\|_-)\frac{a_{-k+j}\varphi_k}{C_-\|f\|_-},
& &
\displaystyle
\sigma_5 :=
\sum_{j=1}^{[k/2]} b_{-j}\|f\|_-, 
\\
\\
\displaystyle
\Sigma_6 :=
\sum_{j<-[k/2]}(b_j\|a\|_-)\cdot 0=0,
& &
\displaystyle
\sigma_6 :=
\sum_{j<-[k/2]} b_j\|f\|_-.
\end{array}
\]
Since $\Phi$ is convex and $\Phi(0)=0$,
from (\ref{eq:main.1}) and Jensen's inequality
(see, e.g., \cite[Theorem~90]{HLP52}) it follows that
\begin{eqnarray}
\label{eq:main.2}
&&
\sum_{k=1}^N
\Phi\left(
\frac{|(fg)_{-k}|\varphi_k}{C_-(\|f\|_W\|g\|_-+\|g\|_W\|f\|_-)}
\right)w_k
\\[3mm]
\nonumber
&&
\le
\frac{\|g\|_-\Big(\sigma_1(N)+\sigma_2(N)\Big)+
\|f\|_-\Big(\sigma_3(N)+\sigma_4(N)\Big)}
{\|f\|_W\|g\|_-+\|g\|_W\|f\|_-},
\end{eqnarray}
where
\begin{eqnarray*}
\sigma_1(N) 
:=
\sum_{k=1}^N \sum_{j=0}^\infty
a_j\Phi\left(\frac{b_{-k-j}\varphi_k}{C_-\|g\|_-}\right)w_k,
&\!\!\!\!\!&
\sigma_2(N)
:=
\sum_{k=1}^N\sum_{j=1}^{[k/2]}
a_{-j}\Phi\left(\frac{b_{-k+j}\varphi_k}{C_-\|g\|_-}\right)w_k,
\\[3mm]
\sigma_3(N) 
:=
\sum_{k=1}^N \sum_{j=0}^\infty
b_j\Phi\left(\frac{a_{-k-j}\varphi_k}{C_-\|f\|_-}\right)w_k,
&\!\!\!\!\!&
\sigma_4(N)
:=
\sum_{k=1}^N\sum_{j=1}^{[k/2]}
b_{-j}\Phi\left(\frac{a_{-k+j}\varphi_k}{C_-\|f\|_-}\right)w_k.
\end{eqnarray*}
Taking into account (i), (ii), and $C_\varphi\ge 1$, 
we have 
\begin{eqnarray}
\label{eq:main.3}
&&
\sigma_1(N) 
=
\sum_{j=0}^\infty a_j \sum_{k=1}^N \Phi\left(\frac{b_{-k-j}\varphi_k}{C_-\|g\|_-}\right)w_k
\le
\sum_{j=0}^\infty a_j \sum_{k=1}^N \Phi\left(\frac{b_{-k-j}\varphi_{k+j}}{C_-\|g\|_-}\right)w_{k+j}
\\[3mm]
\nonumber
&&\le
\sum_{j=0}^\infty a_j \sum_{k=1}^\infty \Phi\left(\frac{b_{-k}\varphi_k}{C_-\|g\|_-}\right)w_k
\le
\|f\|_W\sum_{k=1}^\infty \Phi\left(\frac{C_\varphi b_{-k}\varphi_k}{C_-\|g\|_-}\right)w_k.
\end{eqnarray}

Changing variables in the second sum in $\sigma_2(N)$ ($r=k-j$), we obtain
\[
\sigma_2(N):=
\sum_{k=1}^N\sum_{j=1}^{[k/2]}a_{-j}\Phi\left(\frac{b_{-k+j}\varphi_k}{C_-\|g\|_-}\right)w_k
=
\sum_{k=1}^N\sum_{r=k-[k/2]}^{k-1}a_{r-k}\Phi\left(\frac{b_{-r}\varphi_k}{C_-\|g\|_-}\right)w_k.
\]
Since $\Phi$ is increasing, it follows from the latter equality and Proposition~\ref{pr:3} that
\begin{eqnarray}
\label{eq:main.4}
&&
\sigma_2(N) \le
C_w\sum_{k=1}^N\sum_{j=k-[k/2]}^{k-1} a_{-k+j}\Phi\left(\frac{C_\varphi b_{-j}\varphi_j}{C_-\|g\|_-}\right)w_j
\\[3mm]
\nonumber
&&\le
C_w\sum_{k=1}^N\sum_{j=1}^\infty a_{-k+j}\Phi\left(\frac{C_\varphi b_{-j}\varphi_j}{C_-\|g\|_-}\right)w_j
=
C_w\sum_{j=1}^\infty\Phi\left(\frac{C_\varphi b_{-j}\varphi_j}{C_-\|g\|_-}\right)w_j\sum_{k=1}^N a_{-k+j}
\\[3mm]
\nonumber
&&\le 
C_w\sum_{j=1}^\infty\Phi\left(\frac{C_\varphi b_{-j}\varphi_j}{C_-\|g\|_-}\right)w_j\sum_{k=-\infty}^{\infty}a_{-k}
=
C_w\|f\|_W \sum_{j=1}^\infty\Phi\left(\frac{C_\varphi b_{-j}\varphi_j}{C_-\|g\|_-}\right)w_j.
\end{eqnarray}
Combining (\ref{eq:main.3})--(\ref{eq:main.4}), we arrive at
\begin{equation}\label{eq:main.5}
\sigma_1(N)+\sigma_2(N)
\le
\|f\|_W(1+C_w)
\sum_{k=1}^\infty
\Phi\left(\frac{C_\varphi b_{-k}\varphi_k}{C_-\|g\|_-}\right)w_k.
\end{equation}
Since $\Phi$ is an Orlicz function, $\Phi(x)/x$ is a non-decreasing function
(see, e.g., \cite[p.~139]{LT77}). Thus,
\[
(1+C_w)\Phi(x)\le \Phi((1+C_w)x),
\quad
x\ge 0.
\]
Applying this inequality to (\ref{eq:main.5}), we obtain for $N\ge 1$,
\begin{equation}\label{eq:main.6}
\sigma_1(N)+\sigma_2(N)
\le
\|f\|_W
\sum_{k=1}^\infty 
\Phi\left(\frac{b_{-k}\varphi_k}{\|g\|_-}\right)w_k
\le
\|f\|_W
\end{equation}
Analogously one can show that for $N\ge 1$,
\begin{equation}\label{eq:main.7}
\sigma_3(N)+\sigma_4(N)\le\|g\|_W.
\end{equation}
Taking into account that $N$ is arbitrary, from (\ref{eq:main.2}),
(\ref{eq:main.6}), and (\ref{eq:main.7}) we get
\[
\sum_{k=1}^\infty\Phi\left(
\frac{|(fg)_{-k}|\varphi_k}{C_-(\|f\|_W\|g\|_-+\|g\|_W\|f\|_-)}
\right)w_k
\le 1.
\]
Therefore,
\begin{equation}\label{eq:main.8}
\|fg\|_-
\le 
C_-(\|f\|_W\|g\|_-+\|g\|_W\|f\|_-)
\le
2C_-\|f\|_{W\cap F}\|g\|_{W\cap F}.
\end{equation}
By using of Proposition~\ref{pr:1}(b), one can similarly prove that
\begin{equation}\label{eq:main.9}
\|fg\|_+
\le 
C_+(\|f\|_W\|g\|_++\|g\|_W\|f\|_+)
\le
2C_+\|f\|_{W\cap F}\|g\|_{W\cap F},
\end{equation}
where $C_+:=(1+C_\varrho)C_\psi$, $C_\psi$ and $C_\varrho$
are the constants in the $\Delta_2$-condition for the
sequences $\psi=\{\psi_n\}_{n=0}^\infty$ and
$\varrho=\{\varrho_n\}_{n=0}^\infty$, respectively.
Combining (\ref{eq:main.8}), (\ref{eq:main.9}), and
\[
\|fg\|_W\le\|f\|_W\|g\|_W\le\|f\|_{W\cap F}\|g\|_{W\cap F},
\]
we arrive at (\ref{eq:main*}).
\end{proof}
\subsection*{Acknowledgment}
I would like to thank  Albrecht B\"ottcher
(Chemnitz Technical University, Germany) 
for useful remarks on an earlier version of this paper.

\end{document}